%% file: jpurcell-slope.tex
\documentclass[11pt]{amsart}

\usepackage{amsfonts, amstext, amsmath, amsthm, amscd, amssymb}
\usepackage{graphicx, color, epsfig}

\newtheorem{theorem}{Theorem}[section]
\newtheorem{corollary}[theorem]{Corollary}
\newtheorem{lemma}[theorem]{Lemma}

\newtheorem{proposition}[theorem]{Proposition}

\theoremstyle{definition}
\newtheorem{define}[theorem]{Definition}
\newtheorem{remark}[theorem]{Remark}

\newcommand{\HH}{{\mathbb{H}}}

\newcommand{\CC}{{\mathbb{C}}}

\title{Slope lengths and generalized augmented links}

\author{Jessica S. Purcell}
\thanks{{The author is supported in part by NSF grant DMS--0704359.}}
\address{Mathematical Institute, 24-29 St Giles', Oxford, OX1 3LB,
England}
\address{Department of Mathematics, Brigham Young
University, Provo, UT 84602}
\email{jpurcell@math.byu.edu}


\begin{document}
\bibliographystyle{hamsplain}

\begin{abstract}
In this paper, we determine geometric information on slope lengths of
a large class of knots in the 3--sphere, based only on diagrammatical
properties of the knots.  In particular, we show such knots have
meridian length strictly less than $4$, and we find infinitely many
families with meridian length approaching $4$ from below.  Finally, we
present an example to show that, in contrast to the case of the
regular augmented link, longitude lengths of these knots cannot be
determined by a function of the number of twist regions alone.
\end{abstract}
\maketitle

% body 

% Section: Intro
\section{Introduction}
\label{sec:intro}
\input{intro}

% Section: Define links
\section{Generalized augmented links}
\label{sec:auglink}
\input{auglink}

% Section: Meridian lengths
\section{Meridians}
\label{sec:cusp-shape}
\input{cusp-shape}

% Section: Families with meridian approaching 4
\section{Families with meridians approaching $4$}
\label{sec:families}
\input{families}

% Section: Longitude is independent of twist number
\section{Longitude length is independent of twist number}
\label{sec:longitude}
\input{longitude}

\bibliography{references}

\end{document}

%% file: intro.tex
Given a knot in the 3--sphere whose complement admits a complete
hyperbolic structure, by Mostow--Prasad rigidity the hyperbolic
structure is unique, and therefore a knot invariant.  However, it
seems to be a difficult problem to determine geometric properties of
the hyperbolic structure based only on a diagram of the knot.

In a recent paper \cite{purcell:aug}, we discussed geometric
properties of \emph{generalized augmented links}.  These allowed us to
determine geometric information on large classes of knots in the
3--sphere, including geometric information on highly non-alternating
knots, based only on a diagram.

In this paper, we extend the discussion of the geometry to information
on slope lengths.  In particular, we show such knots have meridian
length strictly less than $4$, and we find infinitely many families
with meridian length approaching $4$ from below.  This is of interest,
because to date there are no known examples of knots with meridian of
length $4$ or more, although Agol has found a family of knots whose
meridian length approaches $4$ from below \cite{agol:bounds}.  This
paper provides evidence that knots with meridian length less than $4$
are extremely common.  Finally, we show that, in contrast to the case
of the regular augmented link, longitude lengths of these knots cannot
be determined by a function of the number of twist regions alone.  We
do this by presenting examples of knots with arbitrarily many twist
regions, but with longitude length bounded by $4$.

The examples of knots of this paper are of interest, as they give
classes of diagrams for which geometric information is available, but
which are highly non-alternating and have arbitrarily large volume.

\subsection{Meridians}

A \emph{slope} is defined to be an isotopy class of simple closed
curves on a torus boundary of a 3--manifold.  When the manifold admits
a hyperbolic structure, the slope inherits a length: the length of a
geodesic representative of the slope in a horospherical torus about
the boundary.

A knot in the 3--sphere has at least one short slope: the meridian of
the knot.  By the 6--Theorem, proved independently by Agol
\cite{agol:bounds} and Lackenby \cite{lackenby:word}, the meridian has
length strictly less than $6$.

Interestingly, to date no knot complements have been found with
meridian length at least $4$.  Agol found a family of examples of
knots with meridian length approaching $4$ from below
\cite{agol:bounds}.  Agol's knots are described as Dehn fillings of
components of the 2--fold branched cover over one component of the
Borromean rings, and in fact, they form a subset of the knots of this
paper, as we will see below.  Adams, and Adams and co-authors have
found upper bounds on meridian lengths of certain classes of knots.
In particular, 2--bridge knots have meridian length less than $2$
\cite{adams:2-gen}, the meridian length of an alternating knot is at
most $3-6/c$, where $c$ is the number of crossings of the knot, and
more generally, the meridian of a knot has length less than $6-7/c$
\cite{adamsetal:cusp-size}.  However, this general bound becomes
little better than that given by the $6$-Theorem as $c$ increases.

In this paper, we provide evidence that $4$ may be a better upper
bound on meridian length than $6$.  We do so by presenting a family of
links, called \emph{generalized augmented links}, from which any knot
may be obtained by Dehn filling.  We show that the slope on these
links which corresponds to the meridian of the knot has length
strictly less than $4$.  Since high Dehn filling yields a manifold
geometrically close to the original (see \cite{thurston}), for any
$\epsilon>0$ we need only exclude a finite number of slopes for each
link component, and the knot obtained by Dehn filling will have
meridian length within $\epsilon$ of that of the generalized augmented
link.  Thus if we enumerate knots by Dehn fillings of generalized
augmented links, ``most'' knots have meridian of length less than $4$.

We call the links \emph{generalized} augmented links because they
generalize a class of links, called augmented links, that has been
studied extensively in the past.  It is well known that any knot can
be obtained by Dehn filling a regular augmented link.  Actually, it
was shown independently by the author with Futer \cite{futer-purcell},
and by Schoenfeld \cite{schoenfeld} that these regular augmented links
have meridians of length exactly $2$.  Note this implies that if we
restrict to augmented links and their Dehn fillings, we obtain an
enumeration of knots for which ``most'' knots actually have meridian
lengths approaching $2$.  

Of course, the result of $2$ is unsatisfactory for many reasons, which
partially motivated this paper.  First and foremost, $2$ is much
shorter than the current largest known meridian length, i.e. the
examples of Agol approaching $4$ \cite{agol:bounds}.  Thus if we wish
to determine an upper bound, or find examples with larger meridians,
we cannot consider high Dehn fillings on regular augmented links.

Secondly, knots obtained by high Dehn fillings of regular augmented
links are very closely related to alternating knots.  Such knots are
either alternating or only non-alternating between twist regions.
That is, they are only non-alternating in very restricted ways.  In
fact, regular augmented links actually give a very nice way of
enumerating alternating knots (see for example \cite{adams:aug},
\cite{lackenby:alt-volume}).  By the above results, any examples of
knots with meridian of length greater than $4$ will most likely come
from classes of highly non-alternating knots.

To that end, we investigate highly non-alternating knots and their
meridian lengths.  We will see that Dehn fillings of generalized
augmented links give a nice description of many highly non-alternating
knots.  In addition, we show that generalized augmented links fit in
with known meridian bounds much better than regular augmented links,
by using them to construct families of knots whose meridian length
approaches $4$ from below.  Thus, to Agol's example of a family of
knots with meridian length approaching $4$ from below, we add
infinitely many more examples of families of knots with meridian
length approaching $4$ from below.

\subsection{Longitudes}

In \cite{futer-purcell}, we showed that the longitude of a regular
augmented link was bounded below by the number of twist regions of a
diagram (we review the definition of twist region in Section
\ref{sec:auglink}).

In this paper, we show that no such result holds for generalized twist
regions and generalized augmented links, by presenting an example of
hyperbolic knots with arbitrarily many twist regions, but longitude
length bounded near $4$.

\subsection{Organization}

This paper is organized as follows.  In Section \ref{sec:auglink} we
define generalized augmented links and present some examples.  In
Section \ref{sec:cusp-shape} we prove that the slope which corresponds
to the meridian has length less than $4$.  In Section
\ref{sec:families} we give examples of families of knots with
meridians approaching length $4$ from below.  Finally, in Section
\ref{sec:longitude} we present examples of knots with arbitrarily many
generalized twist regions, but longitude length bounded near $4$.

%% file: auglink.tex
\subsection{Twisting}

In \cite{adams:aug}, Adams defined a notion of augmented alternating
links, and proved that these were hyperbolic.  We will generalize his
definition.  We also generalize a related notion, that of
twist region of a knot or link diagram.

\begin{figure}
\begin{center}
	\begin{tabular}{ccccc}
			(a) &
	\includegraphics{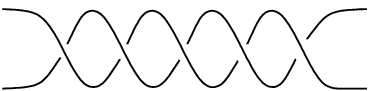} &
	\hspace{.2in} &
	(b) &
	\includegraphics{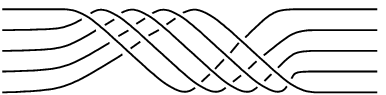}
	\end{tabular}
\end{center}
	\caption{(a) A twist region.  (b) A generalized twist region.}
	\label{fig:twist}
\end{figure}

\begin{define}
	A \emph{twist region} of a diagram is a portion of a
	diagram in which two strands twist about each other maximally, as in
	Figure \ref{fig:twist} (a).

	More formally, a diagram of a knot is a 4--valent graph with
	over--under crossing information at each vertex.  This graph divides
	the projection plane into complementary regions with some number of
	sides.  A \emph{bigon} is a region with just two sides.  A
	\emph{twist region} of a diagram is the portion of a diagram inside
	a simple closed curve which bounds a string of bigons arranged end
	to end.  We require twist regions to be alternating, else we can
	reduce the number of crossings in the region, and we require them to
	be maximal, in the sense that there are no other bigons adjacent to
	either end of the string of bigons in a twist region.  A single
	crossing adjacent to no bigons is also a twist region.
\end{define}

\begin{define}
	A \emph{generalized twist region} of a diagram is a portion of the
	diagram in which $m$ strands twist about each other maximally, for
	$m\geq 2$.  See Figure \ref{fig:twist} (b).

	More precisely, the two strands of a twist region bound an embedded
	\emph{ribbon surface}.  In a generalized twist region, additional
	strands are allowed to run parallel to the two strands, embedded on
	the ribbon surface.  
\end{define}

In a regular twist region on two strands, a \emph{full twist} consists
of two crossings.  Each strand twists once around the other and then
exits the full twist in the same position in which it entered.  A
\emph{half twist} is a single crossing of the two strands.  We
generalize these definitions to generalized twist regions.

\begin{figure}
  \begin{center}
		\includegraphics{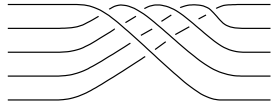}
  \end{center}
	\caption{Five strands make a single half twist.}
\label{fig:full-half}
\end{figure}

\begin{define}
	A \emph{full twist} in a generalized twist region is a region where
	the ribbon surface containing the strands of the twist region makes
	a full twist.  That is, the two outermost strands of the generalized
	twist region form a regular twist region.  In a full twist those two
	strands twist around each other once.  Figure \ref{fig:twist}(b)
	shows a single full twist of five strands.

	Similarly, a \emph{half twist} is a region where the ribbon surface
	of a generalized twist region makes a single half twist.  This
	corresponds to a single crossing of the outermost strands.  See
	Figure \ref{fig:full-half}.
\end{define}

\begin{figure}
\begin{center}
\input{figures/mult-tw-region.pstex_t}
\end{center}
\caption{Five strands make a full and a half twist of all five
strands, followed by a half twist of three strands, and a half twist
of two strands (in the opposite direction).}
\label{fig:twist-ex}
\end{figure}
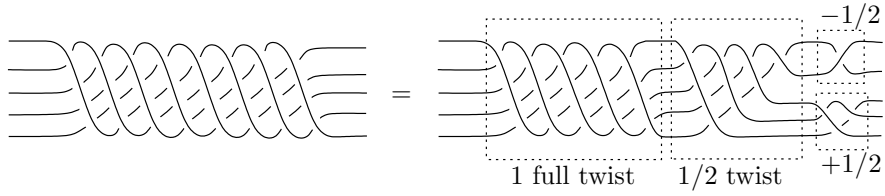

In Figure \ref{fig:twist-ex}, an example is shown of a portion of a
diagram containing three twist regions.  Notice it is highly
non-alternating, and each twist region containing more than $2$
strands is highly non-alternating.  Do note that we could have
selected each single crossing of this portion of the diagram to be a
twist region on two strands.  The choice of twist regions is not
unique.  However, that will not affect our results.

\subsection{Augmenting}

\begin{figure}
\centerline{\includegraphics{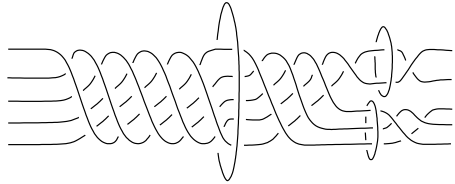}}
\caption{Insert crossing circles into knot diagram at each twist region.}
\label{fig:cross-circ}
\end{figure}

We \emph{augment} a diagram of a knot or link by inserting a simple
closed curve encircling each generalized twist region.  This is called
a \emph{crossing circle}.  See Figure \ref{fig:cross-circ}.  The
components of the link coming from the original link components are
called the \emph{knotting strands}.

When each twist region has just two strands, as in the classical
setting, the resulting link is called an \emph{augmented link}.
Provided a diagram is sufficiently reduced in obvious ways, an
augmented link with at least two twist regions will be hyperbolic with
very explicit geometry.  See for example the papers of Futer and
Purcell \cite{futer-purcell}, and Purcell \cite{purcell:cusps}.

\begin{define}
	A \emph{generalized augmented link} is a link which is formed from a
	knot or link diagram by augmenting each generalized twist region.
\end{define}

Generalized augmented links do not have quite as nice geometry as
regular augmented links, but we can still deduce facts about their
geometry.

\begin{figure}
\begin{center}
  \includegraphics{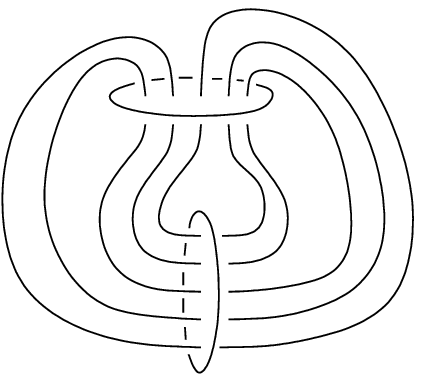}
  \hspace{.25in}
  \includegraphics{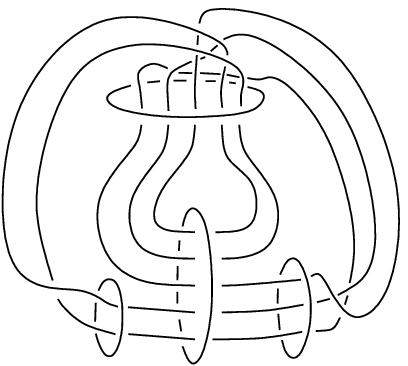}
\end{center}
  \caption{Untwisted generalized augmented links.}
 \label{fig:gen_aug_link}
\end{figure}

Suppose in a knot or link diagram, the $i$-th twist region consists of
$t_i$ full twists, plus possibly a half twist.  We form the
generalized augmented link $L'$ by adding crossing circles.  In
particular, in the $i$-th twist region of the diagram of $L'$, there
are $t_i$ full twists.  Obtain a new link $L$ by removing all $t_i$
full twists from the diagram of $L'$, for each $i$.  Thus $L$ has a
diagram consisting of crossing circle components bounding the
component from the knot.  This knot component is either embedded in
the projection plane, if all the $t_i$ happened to be even, or it may
contain single half twists at crossing circles.  The links of Figure
\ref{fig:gen_aug_link} are of this form.

Note $S^3-L'$ and $S^3-L$ are homeomorphic as manifolds.  This is
because twisting and untwisting along a crossing circle $C_i$ give
homeomorphisms of the solid torus $S^3-C_i$.  If one of $S^3-L$ or
$S^3-L'$ is hyperbolic, Mostow--Prasad rigidity implies that the two
manifolds are isometric.  Thus we will assume from now on that a
generalized augmented link has a diagram with all full twists removed
from each generalized twist region.

Note next that we may obtain the original knot complement by Dehn
filling $S^3-L$ along crossing circles.  In particular, if $t_i$ full
twists were removed from the $i$-th twist region to go from the
diagram of $L'$ to the diagram of $L$, then by performing $1/t_i$ Dehn
filling on the $i$-th crossing circle, for each $i$, we obtain again
the original knot complement.  Consequently, any knot is obtained by
Dehn filling a generalized augmented link.

%% file: figures/mult-tw-region.pstex_t
\begin{picture}(0,0)%
\includegraphics{figures/mult-tw-region.pstex}%
\end{picture}%
\setlength{\unitlength}{3947sp}%
\begingroup\makeatletter\ifx\SetFigFont\undefined%
\gdef\SetFigFont#1#2#3#4#5{%
  \reset@font\fontsize{#1}{#2pt}%
  \fontfamily{#3}\fontseries{#4}\fontshape{#5}%
  \selectfont}%
\fi\endgroup%
\begin{picture}(5563,1162)(1559,-1301)
\put(6676,-1186){\makebox(0,0)[lb]{\smash{{\SetFigFont{10}{12.0}{\familydefault}{\mddefault}{\updefault}{\color[rgb]{0,0,0}$+1/2$}%
}}}}
\put(3976,-736){\makebox(0,0)[lb]{\smash{{\SetFigFont{10}{12.0}{\familydefault}{\mddefault}{\updefault}{\color[rgb]{0,0,0}$=$}%
}}}}
\put(4726,-1261){\makebox(0,0)[lb]{\smash{{\SetFigFont{10}{12.0}{\familydefault}{\mddefault}{\updefault}{\color[rgb]{0,0,0}$1$ full twist}%
}}}}
\put(5776,-1261){\makebox(0,0)[lb]{\smash{{\SetFigFont{10}{12.0}{\familydefault}{\mddefault}{\updefault}{\color[rgb]{0,0,0}$1/2$ twist}%
}}}}
\put(6670,-247){\makebox(0,0)[lb]{\smash{{\SetFigFont{10}{12.0}{\familydefault}{\mddefault}{\updefault}{\color[rgb]{0,0,0}$-1/2$}%
}}}}
\end{picture}%

%% file: cusp-shape.tex
\subsection{Reflection and totally geodesic surfaces}

The key to the meridian slope length estimates for a hyperbolic
generalized augmented link is the existence of embedded totally
geodesic surfaces in the link complement.

\begin{proposition}
  Any generalized augmented link complement admits an involution
	fixing a (possibly disconnected) embedded surface pointwise.
  \label{prop:refl_plane}
\end{proposition}

\begin{proof}
First consider the case where the generalized augmented link has no
half twists.  That is, all link components that are not crossing
circles are embedded on the projection plane, as in the figure on the
left in Figure \ref{fig:gen_aug_link}.  Then there is an involution of
the link complement fixing the projection plane pointwise, obtained by
reflecting the projection plane.  

When there are half twists, the picture is only slightly more
complicated.  In this case, when we reflect in the projection plane,
we obtain a new link in which all crossings of the half twists have
been reversed.  That is, they remain half twists, only now in the
opposite direction of those in the original diagram.  There is a
homeomorphism taking the original link complement to that with half
twists reversed: simply twist one full time around each crossing
circle which contains half twists.  Following the reflection with the
inverse of this homeomorphism gives the desired involution.
\end{proof}

\begin{define}
  The \emph{reflection surface} refers to the surface fixed pointwise by
  the involution of Proposition \ref{prop:refl_plane}.
\end{define}

When the generalized augmented link is hyperbolic, it is a well--known
consequence of the proof of Mostow's theorem that the involution must
be an isometry, and the reflection surface must be totally geodesic.
From now on, we restrict to hyperbolic augmented links, and we
consider the reflection surfaces as totally geodesic embedded surfaces
in the link complement.

Note from the proof of Proposition \ref{prop:refl_plane} that the
reflection surface is just the projection plane when the augmented
link has no half twists.  When the augmented link has half twists, the
reflection surface still corresponds to the projection plane outside a
neighborhood of the disks bounded by the crossing circles.  We call
these disks \emph{crossing disks}.  Notice also that the components of
a generalized augmented link that come from the original (unaugmented)
link components are embedded in the reflection surface.

\begin{lemma}
	Let $K$ be a knotting strand of a hyperbolic generalized
	augmented link.  Then a meridian of $K$ is isotopic to a curve
	perpendicular to the reflection surface, which intersects it twice.
\label{lemma:merid-perp}
\end{lemma}

\begin{proof}
We may isotope the meridian to lie outside all neighborhoods of
crossing disks.  Away from these neighborhoods, the reflection surface
is the projection plane.  The projection plane divides $S^3$ into two
balls, $S^+$ and $S^-$.  In a small neighborhood of a point on $K$, a
meridian runs from the projection plane, up across the top of a
tubular neighborhood of $K$ through $S^+$, then back to the projection
plane.  It then runs along the bottom of a tubular neighborhood of $K$
through $S^-$, closing again at the projection plane.  Thus it meets
the projection plane, hence the reflection surface, twice.

Note we can arrange the meridian to be taken to itself by the
involution of $S^3-L$, with reversed orientation.  Hence it must be
orthogonal to the reflection surface.
\end{proof}

\subsection{The universal cover}

To find bounds on the lengths of meridians, we will need to use the
hyperbolic structure on $S^3-L$.  Let $K$ denote the knotting strand
of $L$.  The length of a meridian will be measured on a \emph{maximal
cusp} about $K$.  Recall that the maximal cusp is obtained by
expanding a horoball neighborhood about the cusp until it meets
itself.

Under the hyperbolic metric, the boundary of a maximal cusp is a
horospherical torus, with an inherited Euclidean metric.  We will find
bounds on the length of the meridian of the torus under this metric.

Consider the upper half space model $\{(w,t) | w\in\CC, t>0\}$ of the
universal cover $\HH^3$ of $S^3-L$.  Conjugate such that a cusp lifts
to the point at infinity.

The reflection surfaces, as totally geodesic surfaces which meet each
cusp, will lift to totally geodesic planes in $\HH^3$ which meet lifts
of cusps.  Thus in particular, there must be lifts of reflection
surfaces to parallel vertical planes in $\HH^3$.

The involution of $S^3-L$ lifts to a reflection of $\HH^3$ through
each vertical plane, whose projection under the covering map is an
isometry.  Thus the Euclidean distance between these lifts of
reflection surfaces is constant.

Now consider the points of $\CC$ on the boundary at infinity of
$\HH^3$ that project to cusps of $S^3-L$.  Lifts of the geodesic
reflection surfaces must also meet at each of these points.  Thus at
each such point we will see totally geodesic planes, which are
Euclidean hemispheres.  Their boundaries on $\CC$ form a collection of
tangent circles and lines.  Because the reflection surfaces are
embedded, the entire collection of lifts of reflection surfaces,
including hemispheres and vertical planes in $\HH^3$, will consist of
disjoint planes, tangent only at points at infinity.

Finally, we add to this picture a collection of horospheres.  A
maximal cusp about the component $K$ of $S^3-L$ will lift to give a
collection of horoballs in $\HH^3$.  The boundaries of these are
horospheres, which are Euclidean spheres, tangent to points of the
boundary at infinity that project to cusps.  Because we are
considering a maximal cusp, the horospheres will be tangent to each
other in pairs.

By Lemma \ref{lemma:merid-perp}, a meridian of $K$ lifts to a curve
which crosses two lifts of reflection surfaces orthogonally.  Thus the
length of a meridian must be exactly twice the Euclidean distance
between lifts of two reflection surfaces, as measured on one of the
horospheres about a point projecting to $K$.  We will find bounds on
this distance.

\begin{theorem}
A meridian of the knotting strand of a hyperbolic generalized
augmented link, as measured on the maximal cusp about that component,
has length strictly less than $4$.
\label{thm:mer-upper}
\end{theorem}

\begin{proof}
Let $L$ be a hyperbolic generalized augmented link.  Consider a
crossing circle $C_1$ of $L$.  This bounds a crossing disk $D$.  The
disk is punctured by the knotting strand $K$ a total of $m$
times, where $m$ is the number of strands of the generalized twist
region, and $m\geq 2$.  Note that a neighborhood of $K$ intersects $D$
in meridians of $K$.  Note also that $D$ is taken to itself under the
involution of $S^3-L$ of Proposition \ref{prop:refl_plane}, and that
$D$ intersects the reflection surface.  Hence the reflection surface
divides $D$ into two disks, an upper half $D^+$ and a lower half
$D^-$.  Similarly, the reflection surface divides each of the meridians
of $K$ on $D$ into two pieces, an upper half and a lower half.  We
will use the disk $D^+$ to bound the length of the upper half of the
meridian.  Since the reflection is an isometry of $S^3-L$, this will
also give a bound on the length of the lower half.

Consider the universal cover $\HH^3$ of $S^3-L$.  Conjugate such that
infinity projects to the cusp corresponding to the crossing circle
$C_1$.  The disk $D^+$ also lifts to the universal cover.  It will
generally not be totally geodesic, but we may isotope the
intersections of $D$ with the reflection surface to be geodesics.  These
are embedded in reflection surfaces, which are totally geodesic, so in
the universal cover these intersections run across embedded totally
geodesic planes.  Since $D^+$ meets $C_1$, the edges of $D^+$ run from
infinity, across a string of tangent geodesic planes which project to
reflection surfaces, and back to infinity.  See Figure \ref{fig:disk}.

\begin{figure}
\centerline{\includegraphics{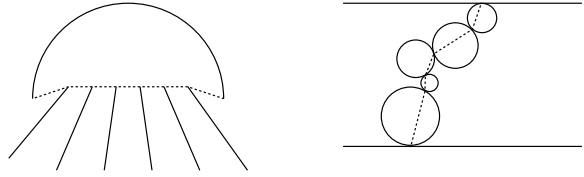}}
\caption{Left: Dotted lines correspond to intersections of $D$ with
	the reflection surface.  Right: In $\HH^3$, these lift to geodesics
	running across lifts of the geodesic reflection surfaces.  Note this
	is a view of $\HH^3$ from infinity.  There are also two dotted lines
	coming out from the horizontal lines to infinity.}
\label{fig:disk}
\end{figure}

Now, as described above, a maximal cusp about $K$ lifts to a
collection of horoballs in $\HH^3$.  We will show that the distance
between two lifts of reflection surfaces, as measured on the
horospherical boundary of one such horoball, is strictly less than
$2$.  Since this length is one half a meridian, it follows that the
length of a meridian is strictly less than $4$.

Our proof is simplified by using an argument of Adams \emph{et al.}
\cite{adams:II}.  In Theorem 4.1 of that paper, they show the width of
a semifree geodesic surface is less than $2$.  That is, they show that
if a manifold contains an embedded totally geodesic surface, if the
surface meets the cusp exactly once, and if there is an embedded disk
in the manifold whose edges lie on the geodesic surface, then the
distance between lifts of the surface on a cusp in $\HH^3$ is strictly
less than $2$.  This is exactly the situation here, except that our
geodesic surface, the reflection surface, meets the cusp more than once.
However, because of the reflective symmetry in our case, the argument
of Adams \emph{et al.}  will still apply.  We reproduce it.

Consider the lift of the disk $D^+$ to $\HH^3$.  In particular,
consider its boundary.  It consists of geodesics embedded on lifts of
reflection surfaces, as in Figure \ref{fig:disk}, alternating with paths
along horospheres at points of tangency of those planes.  Thus we
obtain a string of horospheres, each projecting to the horospherical
torus about the maximal cusp of $K$.  Let $H$ be the horosphere in
this string of smallest Euclidean diameter.  Consider the equator of
$H$, that is, the circle on $H$ of largest diameter parallel to $\CC$.
We will show that the two hemispheres adjacent to $H$, which project
to the reflection surfaces of $S^3-L$, intersect $H$ above its equator.

Suppose not.  Suppose one of the hemispheres, say $S$, intersects $H$
below its equator.  Then $S$ has boundary on $\CC$ a circle of
diameter strictly less than that of $H$.  But then consider the edge
of $D^+$ which lies on $S$.  This runs from $H$ at one end to another
horosphere $H_2$ at the other end.  Because the diameter of $S$ is
less than that of $H$, the diameter of $H_2$ must be less than that of
$H$.  See Figure \ref{fig:triang}.  This contradicts the fact that $H$
was the smallest.

\begin{figure}
\begin{center}
	\input{figures/horo_sizes.pstex_t}
\end{center}
\caption{If $S$ intersects $H$ below its equator, then the radius $r$
	of $S$ is strictly less than the radius $R$ of $H$.  This implies
	that the radius $\widetilde{R}$ of $H_2$ must be strictly less than
	the radius $R$ of $H$.}
\label{fig:triang}
\end{figure}
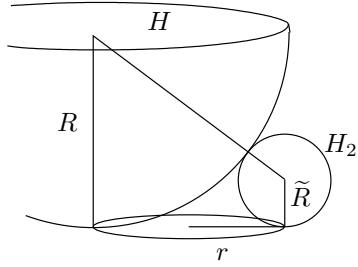

Now, an arc along $H$ running from one lift of a reflection surface
to another lies completely in the upper hemisphere of $H$.  Hyperbolic
geometry shows that the length of an arc running from a point on the
equator to the opposite point on the equator has length exactly $2$.
Thus the distance between reflection surfaces is at most $2$.

Suppose the distance is exactly $2$.  Consider again the string of
horospheres and geodesics bounding the lift of $D^+$ in the universal
cover.  The geodesics must meet the horospheres in their equators,
which implies that each horosphere is exactly the same height and all
the horospheres are tangent.  But then consider $B$, one of the
horospheres lying under a geodesic of the sequence which runs to
infinity.  $\partial D^+ \cap B$ runs between two reflection surfaces,
so has length $2$.  But $\partial D^+ \cap B$ is an arc which must run
from the top of $B$ to its equator.  This has length $1$, giving a
contradiction.

Thus the distance between reflection surfaces is strictly less than $2$.
A meridian runs between two reflection surfaces, so its length is
strictly less than $4$.
\end{proof}

%% file: figures/horo_sizes.pstex_t
\begin{picture}(0,0)%
\includegraphics{figures/horo_sizes.pstex}%
\end{picture}%
\setlength{\unitlength}{3947sp}%
\begingroup\makeatletter\ifx\SetFigFont\undefined%
\gdef\SetFigFont#1#2#3#4#5{%
  \reset@font\fontsize{#1}{#2pt}%
  \fontfamily{#3}\fontseries{#4}\fontshape{#5}%
  \selectfont}%
\fi\endgroup%
\begin{picture}(2430,1670)(350,-1809)
\put(1673,-1769){\makebox(0,0)[lb]{\smash{{\SetFigFont{10}{12.0}{\familydefault}{\mddefault}{\updefault}{\color[rgb]{0,0,0}$r$}%
}}}}
\put(1241,-322){\makebox(0,0)[lb]{\smash{{\SetFigFont{10}{12.0}{\familydefault}{\mddefault}{\updefault}{\color[rgb]{0,0,0}$H$}%
}}}}
\put(676,-961){\makebox(0,0)[lb]{\smash{{\SetFigFont{10}{12.0}{\familydefault}{\mddefault}{\updefault}{\color[rgb]{0,0,0}$R$}%
}}}}
\put(2350,-1089){\makebox(0,0)[lb]{\smash{{\SetFigFont{10}{12.0}{\familydefault}{\mddefault}{\updefault}{\color[rgb]{0,0,0}$H_2$}%
}}}}
\put(2136,-1433){\makebox(0,0)[lb]{\smash{{\SetFigFont{10}{12.0}{\familydefault}{\mddefault}{\updefault}{\color[rgb]{0,0,0}$\widetilde{R}$}%
}}}}
\end{picture}%

%% file: families.tex
In this section we find generalized augmented links whose
meridian length approaches $4$ from below, showing that Theorem
\ref{thm:mer-upper} is sharp.

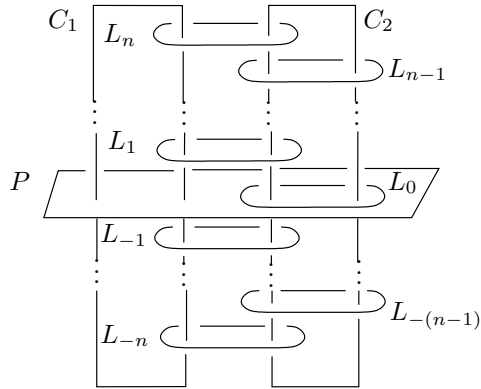
\begin{figure}
\begin{center}
	\input{figures/2bridge-link.pstex_t}
\end{center}
\caption{Dehn filling along $L_i$ and $L_{-i}$, with slope $1/r_i$ and
	$-1/r_i$, respectively, for $i=1, \dots, n$, gives generalized
	augmented link.}
\label{fig:2-bridge}
\end{figure}

The links we consider are links with two crossing circles.  We define
them as Dehn fillings of the link shown in Figure \ref{fig:2-bridge},
as described below.  These links seem to have very interesting
geometric properties in addition to the meridian lengths described
here.  In a forthcoming paper with Futer and Kalfagianni, we will
analyze their volumes \cite{fkp:forthcoming}.

The link in Figure \ref{fig:2-bridge} is a fully augmented 2--bridge
link, with an extra reflective symmetry through the plane $P$ on which
$L_0$ is embedded.  $L_0$ is the link component which will become our
knot.

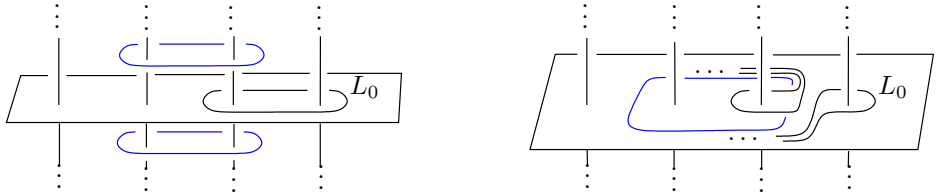
\begin{figure}
\begin{center}
	\input{figures/annular-twist-1.pstex_t}
	\hspace{.5in}
	\input{figures/annular-twist-2.pstex_t}
\end{center}
\caption{The upper and lower link components (left) bound an annulus.
	Dehn filling the top component along slope $1/r_i$, and filling the
	bottom along slope $-1/r_i$ yield $r_i$ Dehn twists of the
	4--punctured sphere (right).}
\label{fig:annular-twist}
\end{figure}

Note that for $1\leq i \leq n$, the link components $L_i$ and $L_{-i}$
bound an annulus in $S^3$.  The Dehn fillings along slope $1/r_i$ on
$L_i$ and $-1/r_i$ on $L_{-i}$ give a twist $r_i$ times about the
annulus.  This restricts to $r_i$ Dehn twists of the horizontal
4--punctured sphere $P$, on which $L_0$ is embedded.  See Figure
\ref{fig:annular-twist}.

Note that if we perform these Dehn fillings in order, first on $L_1$
and $L_{-1}$, then $L_2$ and $L_{-2}$, etc, the result is a link
consisting of two crossing circles, $C_1$ and $C_2$, and a curve $K$
on the plane $P$, which is obtained by a sequence of Dehn twists
applied to $L_0$.  This is a generalized augmented link.  In fact, the
link on the left in Figure \ref{fig:gen_aug_link} has this form.

The geometry of fully augmented 2--bridge links is very well
understood.  In particular, we showed in \cite{purcell:cusps} that a
fully augmented 2--bridge link decomposes into two identical, ideal,
totally geodesic, right angled polyhedra.  We will briefly
explain the decomposition here.  We refer the reader to
\cite[Section 6.4]{purcell:cusps} for more details.

\begin{figure}
\begin{center}
\begin{tabular}{cccc}
	(a) &
	\input{figures/2bridge-decomp2.pstex_t} &
 (b) &
	\input{figures/2bridge-decomp3.pstex_t} 
\end{tabular}
\end{center}
\caption{(a) Cut along the projection plane, and along disks bounded
	by crossing circles.  (b) Retract neighborhood of the link to ideal
	vertices.}
\label{fig:decomp}
\end{figure}
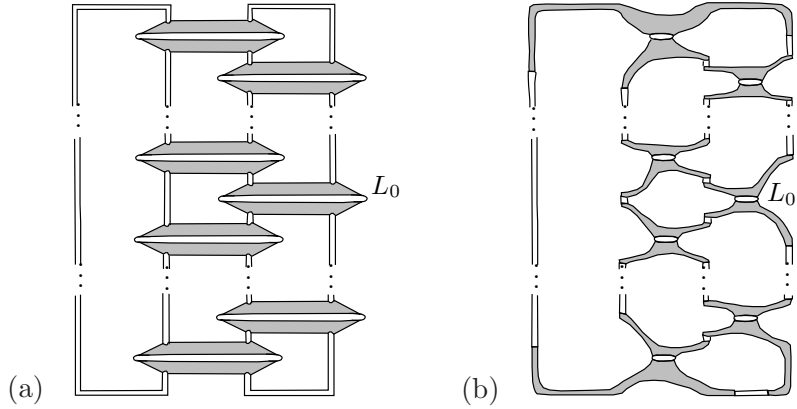

The two polyhedra are obtained by cutting the fully augmented
2--bridge link along the geodesic surface of the projection plane of
Figure \ref{fig:2-bridge}.  This divides the 3--sphere into two
pieces, one above and one below the projection plane.  It divides each
of the 2--punctured disks bounded by crossing circles into two pieces.
Cut along each of these pieces.  See Figure \ref{fig:decomp}(a).

This divides the manifold into two ideal polyhedra.  The ideal
vertices correspond to portions of the link.  See Figure
\ref{fig:decomp}(b).  The faces come from the projection plane as well
as from the 2--punctured disks.  The faces are all totally geodesic,
since the corresponding surfaces are geodesic in the augmented link.
Thus when we give the polyhedra its hyperbolic structure, faces
correspond to (Euclidean) hemispheres in $\HH^3$ orthogonal to the
sphere at infinity.

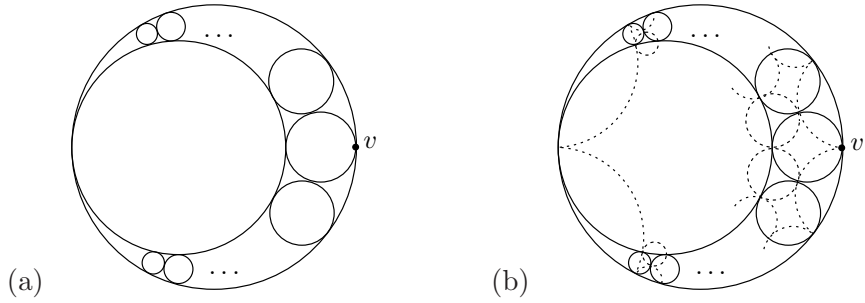
\begin{figure}
\begin{center}
\begin{tabular}{ccccc}
	(a) &
	\input{figures/2bridge-circ1.pstex_t} &
	\hspace{.2in} &
	(b) &
	\input{figures/2bridge-circ2.pstex_t}
\end{tabular}
\end{center}
\caption{Circles giving polyhedral decomposition of augmented
2--bridge link.}
\label{fig:circ}
\end{figure}

One of the two polyhedra is given by the circle packing of Figure
\ref{fig:circ}(a).  The circles should be viewed as boundaries of
spheres on the plane at infinity of $\HH^3$.  They determine totally
geodesic planes in $\HH^3$.  These are the geodesic planes from the
projection plane.  To these circles should be added additional
circles, which meet the given circles exactly in points of tangency,
as in Figure \ref{fig:circ}(b).  These are from the 2--punctured
disks.  Together, these planes cut out an ideal geodesic polyhedron.
To obtain the fully augmented link complement, we may glue two of
these polyhedra together by gluing faces in a manner that reverses the
cutting procedure described above.

We are interested in the cusp shape of the link component $L_0$, since
after high Dehn filling, it will become the knotting strand of our
generalized augmented link (with $C_1$ and $C_2$ crossing circles).  

The vertex marked $v$ in the figure is the only ideal vertex of the
polyhedron that projects to the cusp $L_0$.  Thus the cusp shape of
$L_0$ is obtained by determining the vertex shape of $v$, that is, the
shape of a horosphere intersected with this ideal vertex of the
polyhedron, and gluing two of these together: one for each of the two
polyhedra.

\begin{lemma}
	The shape of the cusp corresponding to $L_0$ is a square.
\label{lemma:square}
\end{lemma}

\begin{figure}
\begin{center}
\begin{tabular}{ccccc}
	(a) &
	\includegraphics{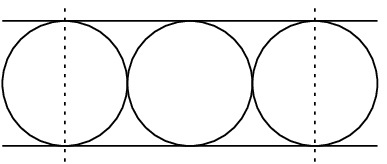}&
	\hspace{.5in}&
	(b)&
	\includegraphics{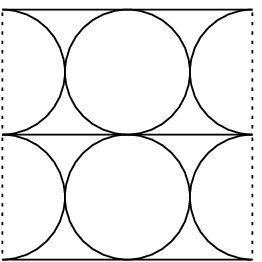}
\end{tabular}
\end{center}
\caption{(a)  The shape of the vertex corresponding to $v$.  (b) The
	shape of the cusp corresponding to $L_0$.}
\label{fig:circ2}
\end{figure}

\begin{proof}
The vertex shape of $v$ is determined by taking the point $v$ of
Figure \ref{fig:circ}(a) to infinity under a M\"obius transformation.
The result is a collection of circles as in Figure \ref{fig:circ2}(a).
Two of these glue together to form the cusp $L_0$, as in Figure
\ref{fig:circ2}(b).  Note that the cusp shape of $L_0$ is a square,
with meridian length $2s$ and longitude $2s$, where $s$ may be
determined by finding the length of a side of the square in a maximal
cusp.
\end{proof}

It remains to determine the maximal cusp.  

\begin{lemma}
Under the metric on a maximal cusp, the meridian and longitude of
$L_0$ both have length $4$.
\label{lemma:L0length}
\end{lemma}

\begin{proof}
Since $v$ is the only vertex of the polyhedron that projects to the
cusp $L_0$, we may expand a horoball about $v$ without it bumping any
of its translates at least until the horoball becomes tangent to a
face of the polyhedron.

Consider again Figure \ref{fig:circ2}(b).  Assume this picture is a
view of $\HH^3$ from infinity.  We may also assume that three corners
of the square lie at $0$, $2$, and $2i$ in $\CC$, so that the circles
of the figure have diameter $1$ in the Euclidean metric on $\CC$.
Thus the corresponding hemispheres have Euclidean radius $1/2$.  Then
a horoball about infinity will not become tangent to one of these
hemisphere faces until its boundary is at height $1/2$.  Hence there
is an embedded horoball of height $1/2$.  If we measure the length of
a meridian and longitude on the corresponding horosphere, the lengths
will be $4$.  This implies that on a maximal cusp the lengths are at
least $4$.

Finally, one may see that this length corresponds to the length on the
maximal cusp by noting that the longitude of $L_0$ bounds a
3--punctured sphere.  It is well known, due to Adams, that the maximal
length of a curve along a 3--punctured sphere is exactly $4$ (see
\cite{adams:3-punct}, or \cite{adams:cusps-collars-systoles}).
\end{proof}

\begin{theorem}
There exist generalized augmented links with meridian lengths
approaching $4$ from below.
\label{thm:mer-upper-sharp}
\end{theorem}

\begin{proof}
For any $\epsilon>0$, and any link of the form in Figure
\ref{fig:2-bridge}, there exist integers $r_1, \dots, r_n$ such that
the generalized augmented link obtained by $1/r_i$ filling on $L_i$
and $-1/r_i$ filling on $L_{-i}$ yields a generalized augmented link
with meridian of length $4-\epsilon$.
\end{proof}

\begin{corollary}
For any $\epsilon>0$, there exist infinitely many knots with two
generalized twist regions whose meridian lengths are greater than
$4-\epsilon$.
\label{cor:length4knots}
\end{corollary}

\begin{proof}
The links of Theorem \ref{thm:mer-upper-sharp} have exactly two
crossing circles $C_1$ and $C_2$.  Find such a link with meridian
length greater than $4-\epsilon/2$.  Then there exist integers $c_1$
and $c_2$ such that Dehn filling $C_1$ along slope $1/c_1$ and Dehn
filling $C_2$ along slope $1/c_2$ yields a knot with meridian length
at least $4-\epsilon$.  
\end{proof}

\begin{remark}
Agol's examples in \cite{agol:bounds} of knots with meridian length
approaching $4$ are actually a subset of the knots of Corollary
\ref{cor:length4knots}.  Agol's knots are obtained by Dehn fillings on
links of the form of Figure \ref{fig:2-bridge} in which $n=1$.  See
\cite[Figure 7]{agol:bounds}.  
\end{remark}

%% file: figures/2bridge-link.pstex_t
\begin{picture}(0,0)%
\includegraphics{figures/2bridge-link.pstex}%
\end{picture}%
\setlength{\unitlength}{3947sp}%
\begingroup\makeatletter\ifx\SetFigFont\undefined%
\gdef\SetFigFont#1#2#3#4#5{%
  \reset@font\fontsize{#1}{#2pt}%
  \fontfamily{#3}\fontseries{#4}\fontshape{#5}%
  \selectfont}%
\fi\endgroup%
\begin{picture}(3234,2421)(2104,-2014)
\put(2347,257){\makebox(0,0)[lb]{\smash{{\SetFigFont{10}{12.0}{\familydefault}{\mddefault}{\updefault}{\color[rgb]{0,0,0}$C_1$}%
}}}}
\put(4479,-52){\makebox(0,0)[lb]{\smash{{\SetFigFont{10}{12.0}{\familydefault}{\mddefault}{\updefault}{\color[rgb]{0,0,0}$L_{n-1}$}%
}}}}
\put(4262,-355){\makebox(0,0)[lb]{\smash{{\SetFigFont{10}{12.0}{\familydefault}{\mddefault}{\updefault}{\color[rgb]{0,0,0}$\vdots$}%
}}}}
\put(3707,-355){\makebox(0,0)[lb]{\smash{{\SetFigFont{10}{12.0}{\familydefault}{\mddefault}{\updefault}{\color[rgb]{0,0,0}$\vdots$}%
}}}}
\put(3177,-354){\makebox(0,0)[lb]{\smash{{\SetFigFont{10}{12.0}{\familydefault}{\mddefault}{\updefault}{\color[rgb]{0,0,0}$\vdots$}%
}}}}
\put(2617,-347){\makebox(0,0)[lb]{\smash{{\SetFigFont{10}{12.0}{\familydefault}{\mddefault}{\updefault}{\color[rgb]{0,0,0}$\vdots$}%
}}}}
\put(4478,-762){\makebox(0,0)[lb]{\smash{{\SetFigFont{10}{12.0}{\familydefault}{\mddefault}{\updefault}{\color[rgb]{0,0,0}$L_0$}%
}}}}
\put(2104,-765){\makebox(0,0)[lb]{\smash{{\SetFigFont{10}{12.0}{\familydefault}{\mddefault}{\updefault}{\color[rgb]{0,0,0}$P$}%
}}}}
\put(2721,-507){\makebox(0,0)[lb]{\smash{{\SetFigFont{10}{12.0}{\familydefault}{\mddefault}{\updefault}{\color[rgb]{0,0,0}$L_1$}%
}}}}
\put(2675,-1091){\makebox(0,0)[lb]{\smash{{\SetFigFont{10}{12.0}{\familydefault}{\mddefault}{\updefault}{\color[rgb]{0,0,0}$L_{-1}$}%
}}}}
\put(2630,-1351){\makebox(0,0)[lb]{\smash{{\SetFigFont{10}{12.0}{\familydefault}{\mddefault}{\updefault}{\color[rgb]{0,0,0}$\vdots$}%
}}}}
\put(3725,-1374){\makebox(0,0)[lb]{\smash{{\SetFigFont{10}{12.0}{\familydefault}{\mddefault}{\updefault}{\color[rgb]{0,0,0}$\vdots$}%
}}}}
\put(4276,-1356){\makebox(0,0)[lb]{\smash{{\SetFigFont{10}{12.0}{\familydefault}{\mddefault}{\updefault}{\color[rgb]{0,0,0}$\vdots$}%
}}}}
\put(3175,-1369){\makebox(0,0)[lb]{\smash{{\SetFigFont{10}{12.0}{\familydefault}{\mddefault}{\updefault}{\color[rgb]{0,0,0}$\vdots$}%
}}}}
\put(4492,-1581){\makebox(0,0)[lb]{\smash{{\SetFigFont{10}{12.0}{\familydefault}{\mddefault}{\updefault}{\color[rgb]{0,0,0}$L_{-(n-1)}$}%
}}}}
\put(2676,-1723){\makebox(0,0)[lb]{\smash{{\SetFigFont{10}{12.0}{\familydefault}{\mddefault}{\updefault}{\color[rgb]{0,0,0}$L_{-n}$}%
}}}}
\put(2695,165){\makebox(0,0)[lb]{\smash{{\SetFigFont{10}{12.0}{\familydefault}{\mddefault}{\updefault}{\color[rgb]{0,0,0}$L_n$}%
}}}}
\put(4326,256){\makebox(0,0)[lb]{\smash{{\SetFigFont{10}{12.0}{\familydefault}{\mddefault}{\updefault}{\color[rgb]{0,0,0}$C_2$}%
}}}}
\end{picture}%

%% file: figures/annular-twist-1.pstex_t
\begin{picture}(0,0)%
\includegraphics{figures/annular-twist-1.pstex}%
\end{picture}%
\setlength{\unitlength}{3947sp}%
\begingroup\makeatletter\ifx\SetFigFont\undefined%
\gdef\SetFigFont#1#2#3#4#5{%
  \reset@font\fontsize{#1}{#2pt}%
  \fontfamily{#3}\fontseries{#4}\fontshape{#5}%
  \selectfont}%
\fi\endgroup%
\begin{picture}(2568,1175)(659,-964)
\put(980,-901){\makebox(0,0)[lb]{\smash{{\SetFigFont{10}{12.0}{\familydefault}{\mddefault}{\updefault}{\color[rgb]{0,0,0}$\vdots$}%
}}}}
\put(2612, 95){\makebox(0,0)[lb]{\smash{{\SetFigFont{10}{12.0}{\familydefault}{\mddefault}{\updefault}{\color[rgb]{0,0,0}$\vdots$}%
}}}}
\put(2057, 95){\makebox(0,0)[lb]{\smash{{\SetFigFont{10}{12.0}{\familydefault}{\mddefault}{\updefault}{\color[rgb]{0,0,0}$\vdots$}%
}}}}
\put(1527, 96){\makebox(0,0)[lb]{\smash{{\SetFigFont{10}{12.0}{\familydefault}{\mddefault}{\updefault}{\color[rgb]{0,0,0}$\vdots$}%
}}}}
\put(967,103){\makebox(0,0)[lb]{\smash{{\SetFigFont{10}{12.0}{\familydefault}{\mddefault}{\updefault}{\color[rgb]{0,0,0}$\vdots$}%
}}}}
\put(2828,-312){\makebox(0,0)[lb]{\smash{{\SetFigFont{10}{12.0}{\familydefault}{\mddefault}{\updefault}{\color[rgb]{0,0,0}$L_0$}%
}}}}
\put(2075,-924){\makebox(0,0)[lb]{\smash{{\SetFigFont{10}{12.0}{\familydefault}{\mddefault}{\updefault}{\color[rgb]{0,0,0}$\vdots$}%
}}}}
\put(2626,-906){\makebox(0,0)[lb]{\smash{{\SetFigFont{10}{12.0}{\familydefault}{\mddefault}{\updefault}{\color[rgb]{0,0,0}$\vdots$}%
}}}}
\put(1525,-919){\makebox(0,0)[lb]{\smash{{\SetFigFont{10}{12.0}{\familydefault}{\mddefault}{\updefault}{\color[rgb]{0,0,0}$\vdots$}%
}}}}
\end{picture}%

%% file: figures/annular-twist-2.pstex_t
\begin{picture}(0,0)%
\includegraphics{figures/annular-twist-2.pstex}%
\end{picture}%
\setlength{\unitlength}{3947sp}%
\begingroup\makeatletter\ifx\SetFigFont\undefined%
\gdef\SetFigFont#1#2#3#4#5{%
  \reset@font\fontsize{#1}{#2pt}%
  \fontfamily{#3}\fontseries{#4}\fontshape{#5}%
  \selectfont}%
\fi\endgroup%
\begin{picture}(2596,1175)(4231,-964)
\put(5274,-176){\makebox(0,0)[lb]{\smash{{\SetFigFont{10}{12.0}{\familydefault}{\mddefault}{\updefault}{\color[rgb]{0,0,0}$\ldots$}%
}}}}
\put(5485,-601){\makebox(0,0)[lb]{\smash{{\SetFigFont{10}{12.0}{\familydefault}{\mddefault}{\updefault}{\color[rgb]{0,0,0}$\ldots$}%
}}}}
\put(6212, 95){\makebox(0,0)[lb]{\smash{{\SetFigFont{10}{12.0}{\familydefault}{\mddefault}{\updefault}{\color[rgb]{0,0,0}$\vdots$}%
}}}}
\put(5657, 95){\makebox(0,0)[lb]{\smash{{\SetFigFont{10}{12.0}{\familydefault}{\mddefault}{\updefault}{\color[rgb]{0,0,0}$\vdots$}%
}}}}
\put(5127, 96){\makebox(0,0)[lb]{\smash{{\SetFigFont{10}{12.0}{\familydefault}{\mddefault}{\updefault}{\color[rgb]{0,0,0}$\vdots$}%
}}}}
\put(4567,103){\makebox(0,0)[lb]{\smash{{\SetFigFont{10}{12.0}{\familydefault}{\mddefault}{\updefault}{\color[rgb]{0,0,0}$\vdots$}%
}}}}
\put(6428,-312){\makebox(0,0)[lb]{\smash{{\SetFigFont{10}{12.0}{\familydefault}{\mddefault}{\updefault}{\color[rgb]{0,0,0}$L_0$}%
}}}}
\put(5675,-924){\makebox(0,0)[lb]{\smash{{\SetFigFont{10}{12.0}{\familydefault}{\mddefault}{\updefault}{\color[rgb]{0,0,0}$\vdots$}%
}}}}
\put(6226,-906){\makebox(0,0)[lb]{\smash{{\SetFigFont{10}{12.0}{\familydefault}{\mddefault}{\updefault}{\color[rgb]{0,0,0}$\vdots$}%
}}}}
\put(5125,-919){\makebox(0,0)[lb]{\smash{{\SetFigFont{10}{12.0}{\familydefault}{\mddefault}{\updefault}{\color[rgb]{0,0,0}$\vdots$}%
}}}}
\put(4580,-901){\makebox(0,0)[lb]{\smash{{\SetFigFont{10}{12.0}{\familydefault}{\mddefault}{\updefault}{\color[rgb]{0,0,0}$\vdots$}%
}}}}
\end{picture}%

%% file: figures/2bridge-decomp2.pstex_t
\begin{picture}(0,0)%
\includegraphics{figures/2bridge-decomp2.pstex}%
\end{picture}%
\setlength{\unitlength}{3947sp}%
\begingroup\makeatletter\ifx\SetFigFont\undefined%
\gdef\SetFigFont#1#2#3#4#5{%
  \reset@font\fontsize{#1}{#2pt}%
  \fontfamily{#3}\fontseries{#4}\fontshape{#5}%
  \selectfont}%
\fi\endgroup%
\begin{picture}(2286,2483)(2591,-2044)
\put(3685,-1365){\makebox(0,0)[lb]{\smash{{\SetFigFont{10}{12.0}{\familydefault}{\mddefault}{\updefault}{\color[rgb]{0,0,0}$\vdots$}%
}}}}
\put(3707,-355){\makebox(0,0)[lb]{\smash{{\SetFigFont{10}{12.0}{\familydefault}{\mddefault}{\updefault}{\color[rgb]{0,0,0}$\vdots$}%
}}}}
\put(3177,-354){\makebox(0,0)[lb]{\smash{{\SetFigFont{10}{12.0}{\familydefault}{\mddefault}{\updefault}{\color[rgb]{0,0,0}$\vdots$}%
}}}}
\put(2617,-347){\makebox(0,0)[lb]{\smash{{\SetFigFont{10}{12.0}{\familydefault}{\mddefault}{\updefault}{\color[rgb]{0,0,0}$\vdots$}%
}}}}
\put(4478,-762){\makebox(0,0)[lb]{\smash{{\SetFigFont{10}{12.0}{\familydefault}{\mddefault}{\updefault}{\color[rgb]{0,0,0}$L_0$}%
}}}}
\put(2630,-1351){\makebox(0,0)[lb]{\smash{{\SetFigFont{10}{12.0}{\familydefault}{\mddefault}{\updefault}{\color[rgb]{0,0,0}$\vdots$}%
}}}}
\put(3175,-1369){\makebox(0,0)[lb]{\smash{{\SetFigFont{10}{12.0}{\familydefault}{\mddefault}{\updefault}{\color[rgb]{0,0,0}$\vdots$}%
}}}}
\put(4204,-355){\makebox(0,0)[lb]{\smash{{\SetFigFont{10}{12.0}{\familydefault}{\mddefault}{\updefault}{\color[rgb]{0,0,0}$\vdots$}%
}}}}
\put(4205,-1356){\makebox(0,0)[lb]{\smash{{\SetFigFont{10}{12.0}{\familydefault}{\mddefault}{\updefault}{\color[rgb]{0,0,0}$\vdots$}%
}}}}
\end{picture}%

%% file: figures/2bridge-decomp3.pstex_t
\begin{picture}(0,0)%
\includegraphics{figures/2bridge-decomp3.pstex}%
\end{picture}%
\setlength{\unitlength}{3947sp}%
\begingroup\makeatletter\ifx\SetFigFont\undefined%
\gdef\SetFigFont#1#2#3#4#5{%
  \reset@font\fontsize{#1}{#2pt}%
  \fontfamily{#3}\fontseries{#4}\fontshape{#5}%
  \selectfont}%
\fi\endgroup%
\begin{picture}(1912,2494)(2583,-2045)
\put(3182,-359){\makebox(0,0)[lb]{\smash{{\SetFigFont{10}{12.0}{\familydefault}{\mddefault}{\updefault}{\color[rgb]{0,0,0}$\vdots$}%
}}}}
\put(3707,-355){\makebox(0,0)[lb]{\smash{{\SetFigFont{10}{12.0}{\familydefault}{\mddefault}{\updefault}{\color[rgb]{0,0,0}$\vdots$}%
}}}}
\put(4096,-811){\makebox(0,0)[lb]{\smash{{\SetFigFont{10}{12.0}{\familydefault}{\mddefault}{\updefault}{\color[rgb]{0,0,0}$L_0$}%
}}}}
\put(4209,-351){\makebox(0,0)[lb]{\smash{{\SetFigFont{10}{12.0}{\familydefault}{\mddefault}{\updefault}{\color[rgb]{0,0,0}$\vdots$}%
}}}}
\put(4201,-1347){\makebox(0,0)[lb]{\smash{{\SetFigFont{10}{12.0}{\familydefault}{\mddefault}{\updefault}{\color[rgb]{0,0,0}$\vdots$}%
}}}}
\put(3676,-1370){\makebox(0,0)[lb]{\smash{{\SetFigFont{10}{12.0}{\familydefault}{\mddefault}{\updefault}{\color[rgb]{0,0,0}$\vdots$}%
}}}}
\put(3162,-1365){\makebox(0,0)[lb]{\smash{{\SetFigFont{10}{12.0}{\familydefault}{\mddefault}{\updefault}{\color[rgb]{0,0,0}$\vdots$}%
}}}}
\put(2626,-1359){\makebox(0,0)[lb]{\smash{{\SetFigFont{10}{12.0}{\familydefault}{\mddefault}{\updefault}{\color[rgb]{0,0,0}$\vdots$}%
}}}}
\put(2608,-369){\makebox(0,0)[lb]{\smash{{\SetFigFont{10}{12.0}{\familydefault}{\mddefault}{\updefault}{\color[rgb]{0,0,0}$\vdots$}%
}}}}
\end{picture}%

%% file: figures/2bridge-circ1.pstex_t
\begin{picture}(0,0)%
\includegraphics{figures/2bridge-circ1.pstex}%
\end{picture}%
\setlength{\unitlength}{3947sp}%
\begingroup\makeatletter\ifx\SetFigFont\undefined%
\gdef\SetFigFont#1#2#3#4#5{%
  \reset@font\fontsize{#1}{#2pt}%
  \fontfamily{#3}\fontseries{#4}\fontshape{#5}%
  \selectfont}%
\fi\endgroup%
\begin{picture}(2071,1802)(2397,-1406)
\put(4238,-512){\makebox(0,0)[lb]{\smash{{\SetFigFont{10}{12.0}{\familydefault}{\mddefault}{\updefault}{\color[rgb]{0,0,0}$v$}%
}}}}
\put(3231,186){\makebox(0,0)[lb]{\smash{{\SetFigFont{10}{12.0}{\familydefault}{\mddefault}{\updefault}{\color[rgb]{0,0,0}$\dots$}%
}}}}
\put(3266,-1306){\makebox(0,0)[lb]{\smash{{\SetFigFont{10}{12.0}{\familydefault}{\mddefault}{\updefault}{\color[rgb]{0,0,0}$\dots$}%
}}}}
\end{picture}%

%% file: figures/2bridge-circ2.pstex_t
\begin{picture}(0,0)%
\includegraphics{figures/2bridge-circ2.pstex}%
\end{picture}%
\setlength{\unitlength}{3947sp}%
\begingroup\makeatletter\ifx\SetFigFont\undefined%
\gdef\SetFigFont#1#2#3#4#5{%
  \reset@font\fontsize{#1}{#2pt}%
  \fontfamily{#3}\fontseries{#4}\fontshape{#5}%
  \selectfont}%
\fi\endgroup%
\begin{picture}(2076,1802)(2392,-1406)
\put(3231,186){\makebox(0,0)[lb]{\smash{{\SetFigFont{10}{12.0}{\familydefault}{\mddefault}{\updefault}{\color[rgb]{0,0,0}$\dots$}%
}}}}
\put(3266,-1306){\makebox(0,0)[lb]{\smash{{\SetFigFont{10}{12.0}{\familydefault}{\mddefault}{\updefault}{\color[rgb]{0,0,0}$\dots$}%
}}}}
\put(4238,-525){\makebox(0,0)[lb]{\smash{{\SetFigFont{10}{12.0}{\familydefault}{\mddefault}{\updefault}{\color[rgb]{0,0,0}$v$}%
}}}}
\end{picture}%

%% file: longitude.tex
With Futer in \cite{futer-purcell}, we showed that the length of a
longitude of a regular augmented link, i.e. one with just two strands
per twist region, is bounded below by a linear function of the number
of crossing circles. 

In this section we give an example to prove that similar simple lower
bounds will not be possible in the case of generalized augmented
links.  In particular, we show there exists a family of generalized
augmented links with arbitrarily many twist regions, but longitude
length at most $4$.

\begin{figure}
\begin{center}
  \includegraphics{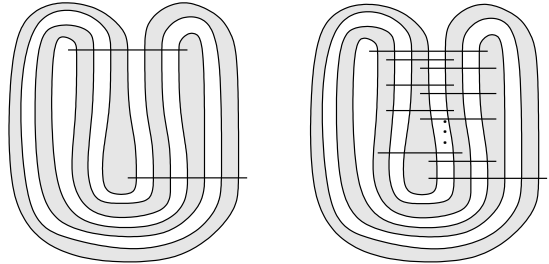}
  \hspace{.2in}
  \input{figures/7strands-ncc.pstex_t}
\end{center}
\caption{Left: The original link $L_2$ with two crossing circles.
	Crossing circles are represented by horizontal lines.
	Right: The generalized augmented link $L_{2n}$ obtained by adding
	$n$ pairs of new crossing circles to $L_2$.}
\label{fig:longitude}
\end{figure}

The example is given by links illustrated in Figure
\ref{fig:longitude}.  In particular, we start with a link $L_2$
consisting of a knotting strand and two crossing circles, $C_1$ and
$C_2$, as shown on the left of that figure.  This link $L_2$ can also
be seen as a 2--bridge link with an extra augmentation component,
exactly as in the previous section, in which the three component links
after Dehn twisting along annuli were 2--bridge links with an extra
augmentation component.  That is, in the link $L_2$, $C_1$ and $C_2$
form the components of a 2--bridge link with reflective symmetry.  The
knotting strand is an augmentation lying flat in the plane of the
reflective symmetry.  Since such 2--bridge links and their
augmentations are known to be hyperbolic, $S^3-L_2$ must be
hyperbolic.  (This is a special case of a theorem on alternating
knots, due originally to Menasco \cite{menasco:alt}.  See also the
paper of Futer and Gueritaud for a more direct proof
\cite{gueritaud-futer}, and Adams for the proof that augmentations are
hyperbolic \cite{adams:aug}.)

Note the knotting strand and crossing circles divide the projection
plane of the link $L_2$ into two distinct 3--punctured spheres.  One
of these 3--punctured spheres is shaded in Figure \ref{fig:longitude}.
To form the new link $L_{2n}$, add $n$ pairs of crossing circles as
shown on the right of Figure \ref{fig:longitude}.  These crossing
circles are chosen so that their intersection with the projection
plane is in the same 3--punctured sphere.  (In the figure, their
endpoints are in the shaded 3--punctured sphere.)  Note the link is a
generalized augmented link, with reflective symmetry through the
projection plane.

\begin{lemma}
The generalized augmented link $L_{2n}$ formed as above is hyperbolic.
\label{lemma:long-hyp}
\end{lemma}

We will prove this lemma at the end of the section.  First, given this
lemma, the results on longitude length are immediate.

\begin{proposition}
	The length of the longitude of the knotting strand component of
	$S^3-L_{2n}$ is at most $4$.
\label{prop:longitude}
\end{proposition}

\begin{proof}
	A longitude bounds a 3--punctured sphere on one side (i.e. the
	unshaded portion of the projection plane in Figure
	\ref{fig:longitude}).  It is well known, due to Adams, that the
	length along a maximal cusp of a hyperbolic 3--punctured sphere is
	at most $4$ (see \cite{adams:3-punct} or
	\cite{adams:cusps-collars-systoles}).  Hence the longitude length is
	at most $4$.
\end{proof}

Proposition \ref{prop:longitude} shows that there can be no
relationship between longitude length and number of generalized twist
regions in a knot diagram.

It remains to show that $S^3-L_{2n}$ is hyperbolic, that is, to prove
Lemma \ref{lemma:long-hyp}.  To do so, we show that $S^3-L_{2n}$ is
irreducible, boundary irreducible, atoroidal, with no essential
annuli.  Hyperbolicity then follows from work of Thurston
\cite{thurston:bulletin}.

Before we begin the proofs, note that we can decompose $S^3-L_{2n}$ in
a simple way.  In particular, if we cut along the projection plane and
the crossing disks (i.e. the disks bounded by the crossing circles),
the result is (topologically) two 3--balls.  Thus these surfaces chop
the manifold $S^3-L_{2n}$, and any submanifolds, into nice pieces.  We
use this fact repeatedly in the lemmas below.

Therefore, to simplify notation, we will denote the projection plane
of $S^3-L_{2n}$ by $P$.  Note it has two components, one of which is a
3--punctured sphere.  The two crossing circles $C_1$ and $C_2$ of
$L_2$ bound crossing disks $E_1$ and $E_2$, respectively.  That is,
$E_1$ is a disk in $S^3$ with boundary the longitude of $C_1$, and
$E_1$ meets the knotting strand of $L_2$ seven times.  Similarly for
$E_2$.  The newly added crossing circles are denoted $C_3, C_4, \dots,
C_{2n}$.  Denote the crossing disk bounded by $C_j$ by $D_j$, for
$j=3, \dots, 2n$.  Each $D_j$ has boundary a longitude of $C_j$, and
meets the knotting strand of $L_{2n}$ four times.  With this notation,
we are ready to begin.

To prove $S^3-L_{2n}$ is irreducible and boundary irreducible, we
first prove the following simpler lemma.

\begin{lemma}
There does not exist a disk in $S^3-L_{2n}$ whose boundary is a
longitude of a crossing circle.
\label{lemma:C_j-non-triv}
\end{lemma}

\begin{proof}
There cannot exist a disk in $S^3-L_{2n}$ whose boundary is a
longitude of $C_1$ or $C_2$, for such a disk would be a boundary
compressing disk of $S^3-L_2$.  Since $S^3-L_2$ is hyperbolic, no such
disk exists.

Suppose there exists a disk $D$ whose boundary is a longitude of
$C_j$, some $j\geq 3$.  Consider the intersections of $D$ with the
disks $E_1$ and $E_2$.  If $D\cap(E_1\cup E_2)$ is non-empty, then
some curve of the intersection bounds a subdisk $E$ in $D$.  We may
assume $E$ meets the projection plane, for otherwise $\partial E$
bounds a subdisk of $E_1$ or $E_2$ which meets no knotting strands,
hence is embedded in $S^3-L_{2n}$.  We may therefore isotope off,
decreasing the number of components of $D\cap (E_1\cup E_2)$.  So
assume $E$ meets $P$.

Now, after isotopy reducing intersections bounding disks, $E\cap P$
consists of an arc running from $E_i$ back to itself.  Since each arc
of the knotting strand runs from $E_1$ to $E_2$, no knotting strand is
bounded between this arc and $E_i$.  So $\partial E$ on $E_i$ bounds a
disk which does not meet any knotting strands.  Thus it bounds a disk
on $E_i$ which is embedded in $S^3-L_{2n}$.  Again we may use this
disk to push the intersection off $E_i$, decreasing the number of
components of $D\cap (E_1\cup E_2)$.

Repeating this process, we may assume $D$ does not meet either $E_i$.
But the intersection of $D$ with the projection plane is then an arc
connecting the two points of intersection of the crossing circle $C_j$
with the projection plane, which does not run through either $C_1$ or
$C_2$.  This is impossible.
\end{proof}

\begin{lemma}
$S^3-L_{2n}$ is irreducible.
\label{lemma:irreducible}
\end{lemma}

\begin{proof}
Suppose $S$ is a 2--sphere in $S^3-L_{2n} \subset S^3-L_2$.  Because
$S^3-L_2$ is irreducible, $S$ must bound a ball $B$ in $S^3-L_2$.  If
$S$ does not bound a ball in $S^3-L_{2n}$, then one of the components
$C_j$ must be contained in $B$.  But $C_3\cup \dots \cup C_{2n}$ is a
trivial link in $S^3$, so it must be homotopically trivial in
$S^3-L_{2n}$.  By the loop theorem, it bounds a disk whose boundary is
a longitude of $\partial N(C_j)$.  This contradicts Lemma
\ref{lemma:C_j-non-triv}.
\end{proof}

\begin{lemma}
	$S^3-L_{2n}$ is boundary irreducible.
\label{lemma:bdry-irred}
\end{lemma}

\begin{proof}
Suppose $D$ is a compressing disk for the boundary of $S^3-L_{2n}$.
By Lemma \ref{lemma:C_j-non-triv}, $\partial D$ is not a longitude of
a crossing circle.  Since we are in $S^3$, $\partial D$ must be
homotopically trivial on $\partial N(L_{2n})$, and we may assume it
bounds a disk $E$ on this boundary surface.  Then $E \cup D$ is a
sphere $S$ in $S^3-L_{2n}$.  By Lemma \ref{lemma:irreducible}, $S$
must bound a ball.  But then we may isotope $D$ through this ball to
lie on $E \subset \partial N(L_{2n})$, and $D$ is not a compressing
disk.
\end{proof}

\begin{lemma}
	$S^3-L_{2n}$ is atoroidal.
\label{lemma:atoroidal}
\end{lemma}

\begin{proof}
Suppose not.  Let $T$ be an essential torus in $S^3-L_{2n}$.  We may
assume first that $T$ meets the projection plane $P$, else it is an
essential torus in one of the two handlebodies of $S^3-L_{2n}-P$.  No
such torus exists.

Then by the equivariant torus theorem, we may assume that $T$ is
preserved under the reflection in the projection plane.  

We obtain a cell decomposition of $T$ by letting edges be
intersections of $T$ with the projection plane, $T\cap P$, and
intersections of $T$ with the crossing disks $T\cap E_1$, $T\cap E_2$,
and $T\cap D_3, \dots, T\cap D_{2n}$.  We let vertices be
intersections of $T$ with $P \cap E_i$, $i=1, 2$, and $P\cap D_j$,
$j=3, \dots, 2n$.  The faces lie in the two balls $S^3-L_{2n}-P-(\cup
E_i)-(\cup D_j)$.

Note each vertex of $T$ is 4--valent.  Thus if we let $v$ denote the
number of vertices, $e$ denote the number of edges, $2e=4v$, or
$e=2v$.  Next note there are no bigon faces or triangle faces of $T$.
Let $f$ denote the number of faces.  Since each edge is adjacent to
two faces, $2e \geq 4f$.  It follows from the Euler characteristic
that every face must be a quadrilateral.

Now, first suppose $T$ meets none of the crossing disks $D_3, D_4,
\dots, D_{2n}$.  $T$ is a torus in $S^3-L_2$, which is hyperbolic, so
$T$ must be compressible or boundary parallel in that manifold.
Consider the intersections of $T$ with the projection plane of
$S^3-L_2$.  Both components of the projection plane are 3--punctured
spheres in $S^3-L_2$, so the possible curves of intersection of $T$
with the projection plane components are limited.  If a curve of
intersection is trivial in the 3--punctured sphere, it remains trivial
when we remove $C_3, \dots, C_{2n}$, since $T$ doesn't meet the
crossing disks bounded by these circles.  Thus in this case $T$ will
be compressible in $S^3-L_{2n}$.  So a curve of intersection must be
non-trivial in the 3--punctured sphere.  Then $T$ is boundary parallel
in $S^3-L_2$.  In fact, it must be parallel to $C_1$ or $C_2$, since
it doesn't meet $D_3, \dots, D_{2n}$.

Since $T$ is essential in $S^3-L_{2n}$, some $C_j$, $j=1, \dots 2n$,
must lie inside the solid torus $V$ bounded by $T$ containing $C_1$ or
$C_2$.  But note $C_1$ and $C_2$ both have endpoints lying in distinct
components of the projection plane, whereas $C_j$, $j=3, \dots, 2n$
have endpoints lying in the same component of the projection plane.
Thus $C_j$ must lie in a ball inside $V$.  This is a ball in $S^3$.
$C_j$ cannot be homotopically trivial, by Lemma
\ref{lemma:C_j-non-triv}, so it must link $C_1$ or $C_2$ inside this
ball.  But this is impossible: the $C_j$, $j=1, \dots, 2n$ are all
unlinked in $S^3$.

So $T$ must meet a crossing disks $D_i$.  Consider $T\cap D_i$.  By
equivariance of $T$, this consists of closed curves encircling some
number of intersections of $D_i$ with the knotting strands.  Since
there are just four such intersections, the curve must encircle $1$,
$2$, $3$, or $4$ of these intersections.

In the case $T\cap D_i$ encircles $1$ knotting strand: Faces adjacent
to the edges $T\cap D_i$ are quadrilaterals.  Two edges must run from
$D_i$ along the projection plane to the next crossing disk, where they
are joined by one other edge.  Each edge along the projection plane
runs parallel to the knotting strand, thus in the next crossing disk
we see that $T$ again encircles $1$ knotting strand.  Continue.  It
follows that $T$ is boundary parallel, parallel to the knotting
strand.

In the case $T\cap D_i$ encircles $3$ knotting strands: Some adjacent
face cannot be a rectangle, for one of the edges running from $D_i$
runs to a crossing disk $D_{i+1}$, say, while the other runs to a
different crossing disk.  Thus the endpoints of these edges cannot be
joined by a single edge.  This is a contradiction.

In the case $T\cap D_i$ encircles $4$ knotting strands, either $T$ is
boundary parallel, parallel to $C_i$, or we run into the same issue as
in the previous case: the two adjacent edges terminate on distinct
crossing disks, so the adjacent face cannot be a quadrilateral.

In the case $T\cap D_i$ encircles $2$ knotting strands, the torus $T$
must enclose the component of $P$ which is a 3--punctured sphere in
$S^3-L_{2n}$.  Adjacent faces to $D_i$ may be quadrilaterals
initially.  However, eventually we will find a face must be adjacent
to the lower crossing disk $C_2$, with one edge running parallel to
the knotting strand to terminate on one side of $E_1$, and the other
running parallel to the knotting strand to terminate on the other side
of $E_1$ in the same component of $P$ minus the knotting strands.
When we connect the endpoints of these edges by an arc, and reflect
through the projection plane, we find they bound a disk on $E_1$.
Hence $T$ is compressible.
\end{proof}

\begin{lemma}
There are no essential annuli in $S^3-L_{2n}$.  
\label{lemma:annulus}
\end{lemma}

\begin{proof}
This follows from Lemma \ref{lemma:atoroidal}, for assume $A$ is an
essential annulus.  If $A$ has boundary components on two distinct
link components, take the boundary of a small embedded neighborhood
of $A$ and of the two link components of $L_{2n}$ on which $\partial
A$ lies.  This is a torus in $S^3-L_{2n}$, hence is inessential.  It
cannot be boundary parallel, so is compressible.  The compressing disk
cannot be on the side of the torus containing $A$ and the two link
components, hence is on the outside.  Cut along this disk.  The result
is a sphere in $S^3-L_{2n}$.  By Lemma \ref{lemma:irreducible}, it
bounds a ball in $S^3-L_{2n}$.  Since the ball can't be on the side of
the sphere containing $A$ and the two link components, it must be on
the other side.  But then the manifold $S^3-L_{2n}$ has just two
boundary components.  This is a contradiction.

If $A$ has boundary components on the same link component, then the
boundary of a small neighborhood of $A$ and the link component gives
two tori in $S^3-L_{2n}$, neither of which may be essential by Lemma
\ref{lemma:atoroidal}.  If either is compressible, consider the
boundary of the compressing disk $D$.  $\partial D$ cannot meet $A$ in
a closed curve, else $A$ is compressible.  Thus $\partial D$ meets $A$
in an arc.  If the arc has both endpoints on the same component of
$\partial A$, then a portion of the disk $D$ and a portion of $A$
together form a disk which must be isotopic to $\partial(S^3-L_{2n})$,
by Lemma \ref{lemma:bdry-irred}, so we may slide this intersection off
of $A$.  Thus we assume the arc of intersection of $\partial D$ and
$A$ is an essential arc on $A$.  Then when we surger along $D$, we
obtain a boundary compressing disk for $S^3-L_{2n}$, hence by Lemma
\ref{lemma:bdry-irred}, the result bounds a ball and so $A$ is
isotopic to $\partial(S^3-L_{2n})$.

Thus when we take the boundary of a small regular neighborhood of $A$
and the link component on which $\partial A$ lies, the two resulting
tori must be boundary parallel.  But then there is some link component
of $L_{2n}$ lying inside.  Then there exists an essential annulus with
boundary components on this link component and on the link component
on which the boundary of $A$ lies.  We are then in the case of the
first paragraph of this proof, and we have a contradiction.
\end{proof}

%% file: figures/7strands-ncc.pstex_t
\begin{picture}(0,0)%
\includegraphics{figures/7strands-ncc.pstex}%
\end{picture}%
\setlength{\unitlength}{3947sp}%
\begingroup\makeatletter\ifx\SetFigFont\undefined%
\gdef\SetFigFont#1#2#3#4#5{%
  \reset@font\fontsize{#1}{#2pt}%
  \fontfamily{#3}\fontseries{#4}\fontshape{#5}%
  \selectfont}%
\fi\endgroup%
\begin{picture}(1510,1643)(68,-966)
\put(901,-211){\makebox(0,0)[lb]{\smash{{\SetFigFont{10}{12.0}{\familydefault}{\mddefault}{\updefault}{\color[rgb]{0,0,0}$\vdots$}%
}}}}
\end{picture}%